\documentclass[12pt]{amsart}
\usepackage{amsmath,amssymb,mathrsfs,hyperref,color}

\usepackage[nameinlink]{cleveref}
\hypersetup{colorlinks={true},linkcolor={blue},citecolor=green}
\crefformat{equation}{(#2#1#3)}

\usepackage{mathtools}
\usepackage[x11names]{xcolor}
\newtagform{blue}{\color{blue}(}{)}

\usepackage{paralist}

\usepackage{subfigure}
\usepackage{float}
\usepackage{palatino}
\usepackage{tikz}
\usepackage{cases}
\usetikzlibrary{intersections}
\usepackage{xcolor}
\usepackage[latin1]{inputenc}
\usepackage{bbm}
\usetikzlibrary{shadings}
\usetikzlibrary[intersections]
\usetikzlibrary{arrows,shapes}

\makeatletter
\def\setliststart#1{\setcounter{\@listctr}{#1}%
  \addtocounter{\@listctr}{-1}}
\makeatother

\makeatletter
\@addtoreset{figure}{section}
\makeatother

\usepackage{calc}
\newtheorem{theorem}{Theorem}[section]
\newtheorem{lemma}[theorem]{Lemma}
\newtheorem{proposition}[theorem]{Proposition}
\newtheorem{corollary}[theorem]{Corollary}
\newtheorem{remarks}[theorem]{Remark}

\newtheorem{definition}[theorem]{Definition}
\numberwithin{equation}{section}

\newcommand{\T}{\mathbb{T}}
\newcommand{\R}{\mathbb{R}}

\newcommand{\N}{\mathbb{N}}

\newcommand{\E}{\mathscr{E}}
\newcommand{\K}{\mathcal{K}}
\newcommand{\PP}{\mathcal{P}}

\newcommand{\M}{\mathcal{M}}

\newcommand{\II}{\mathscr{I}}

\DeclareMathOperator*{\supp}{spt}

\DeclareMathOperator*{\ddiv}{div}

\DeclareMathOperator*{\eps}{\varepsilon}

\makeatletter
\def\moverlay{\mathpalette\mov@rlay}
\def\mov@rlay#1#2{\leavevmode\vtop{%
   \baselineskip\z@skip \lineskiplimit-\maxdimen
   \ialign{\hfil$\m@th#1##$\hfil\cr#2\crcr}}}
\newcommand{\charfusion}[3][\mathord]{
    #1{\ifx#1\mathop\vphantom{#2}\fi
        \mathpalette\mov@rlay{#2\cr#3}
      }
    \ifx#1\mathop\expandafter\displaylimits\fi}
\makeatother

\title[Nash equilibria and Mather measures]{Differential N-players game: Nash equilibria and Mather measures}
\author{Cristian Mendico}
\address{Dipartimento di matematica, Universit\'a degli studi di Roma Tor Vergata -- Via della Ricerca Scientifica 1, 00133 Roma}
\email{mendico@axp.mat.uniroma2.it}

\date{\today}
\thanks{{\bf Acknowledgement:} Cristian Mendico was partly supported by Istituto Nazionale di Alta Matematica (GNAMPA 2022 Research Projects).}

\begin{document}
\usetagform{blue}

\maketitle

\begin{abstract}
We study Nash equilibria for the deterministic ergodic N-players game. We introduce pure strategies, mixed strategies and Nash equilibria associated with those. We show that a Nash equilibrium in mixed strategies exists and it is a Mather measure for the Lagrangian system defined by the cost functional. In conclusion, we show that the mean field limit of the N-players game is described by the ergodic PDE's system for a continuum of players. 
\end{abstract}

\section{Introduction}

Mean field game (MFG) theory, introduced in \cite{bib:LL1, bib:LL2, bib:LL3} and independently in \cite{bib:HCM2, bib:HCM1}, is devoted to the study of differential games with a very large number of interacting agents. A typical model is described by a system of PDEs: an Hamilton-Jacobi equation that provides the optimal choices of players and a Fokker-Planck equation which describes how the distribution of individuals changes. The two equations are coupled by a function that takes into account both the state of a single agent and how he/she is influenced by the others. For an extensive and detailed introduction to the subject we refer to  \cite{bib:BFY}, \cite{bib:NC}, \cite{bib:CFD1, bib:CFD2}, \cite{bib:DEV} and references therein. 

The PDEs system describes a model with a continuum of players which is, clearly, not realistic. However, the solution of such a system is expected to capture the behavior of Nash equilibria for differential N-players game as the number of agents goes to infinity. This problem is generally open and it has been solved in some specific framework, we refer for instance to \cite{bib:ABTHL, bib:CMP, bib:MFJS, bib:DLL, bib:EF, bib:MBEF}. In particular, these works mainly concern with second-order MFG, i.e., systems in which the dynamic of players is described by a controlled stochastic differential equation.  More recently, such a problem have been addressed by analyzing the so-called master equation, that is, an Hamilton-Jacobi equation defined on the Wasserstein space of probability measures. For a comprehensive introduction to such equation and the study of the approximation of Nash equilibria via solutions to the master equation we refer to \cite{bib:CPDLL}.  As far as the authors know, no results are available for first-order MFG.

In this work we address the analysis of first-order differential N-players game with ergodic cost, i.e., given $N \in \N$ individuals we consider the cost
\begin{equation*}
\liminf_{T \to \infty} \frac{1}{T}\int_{0}^{T}{\big(L^{i}(\gamma_{i}(t), u_{i}(t)) + F^{i}(\gamma_{1}(t), \dots, \gamma_{N}(t)) \big)\ {\rm dt}}
\end{equation*}
with dynamic $\dot\gamma_{i}(t)=u_{i}$. For this game we introduce pure strategies and we study the existence of Nash equilibria associated with such strategies. Specifically, we show that if a suitable system of Hamilton-Jacobi equations has a smooth solution then there exists at least one Nash equilibrium in pure strategies and the value of the game is reached at the critical constant of the Lagrangian system defined by the cost. 

However, in general, the system of Hamilton-Jacobi equations mentioned before fails to have a smooth solution and thus the use of pure strategies to find a Nash equilibrium for the game is not satisfactory. Hence, inspired by weak KAM theory \cite{bib:FA}, we define mixed strategies for the game as invariant probability measures for the Euler flow associated with the cost. By using such strategies we find the existence of Nash equilibria under general assumptions and, moreover, we show that Nash equilibria in mixed strategies are Mather measures. 

Finally, considering symmetric games we address the problem of convergence of the N-players game to a suitable PDEs system as N goes to infinity. We end up showing that the limit system is given by the ergodic MFG system already introduced in \cite{bib:CAR, bib:CCMW, bib:CCMW1}. In particular, such a system capture the long-time average behavior of solutions to the classical MFG. 

The paper is organized as follows. Section \ref{sec:weakKAM} is devoted to a preliminary introduction to the classical weak KAM theory for Tonelli Lagrangian (Hamiltonian) systems. In Section \ref{sec:Purestrategies} we introduce the main assumptions on the game, we define pure strategies and Nash equilibria associated with such strategies. We also prove that if a certain Hamilton-Jacobi equation has a smooth solution then there exists at least one Nash equilibrium for the game. In Section \ref{sec:Mathertheory} we show the relation between the ergodic game and Mather theory for Tonelli Hamiltonian systems. Having at our disposal such a relation, in Section \ref{sec:Mixedstrategies} we introduce Nash equilibria in mixed strategies and we show that there always exists a Nash equilibrium and, moreover, we prove that such equilibrium is a Mather measure for the underlying dynamical system. Finally, in Section \ref{sec:Meanfieldlimit} we show that as the number of players goes to infinity the solution to the N-players game converges to a solution to a suitable PDE system of MFG type. We conclude the paper stating a possible future direction of this work towards a new understand of the Nash equilibria for N-player games in Section \ref{sec:Open}.

\section{Known facts on weak KAM theory}
\label{sec:weakKAM}

\begin{definition}[Tonelli Lagrangians]\label{def:def1}
	A function $L: \T^{d} \times \mathbb{R}^{d} \to \mathbb{R}$ is called a {\it Tonelli Lagrangian} if it belongs to $C^{2}$ and it satisfies the following. 
	\begin{itemize}
		\item[(i)] For each $(x,v) \in \T^{d} \times \mathbb{R}^{d} $, the Hessian $D^{2}_{vv}L(x,v)$ is positive definite.
		\item[(ii)] For each $A>0$ there exists $B(A) \in \mathbb{R}$ such that 
		\begin{equation*}
			L(x,v) >A|v|+B(A), \quad  \forall (x,v) \in \T^{d} \times \mathbb{R}^{d}.
		\end{equation*}
		\item[(iii)] For each $R >0$ 
		\begin{equation*}
			A(R):=\sup\Big\{L(x,v): |v| \leq R\Big\} < +\infty.
		\end{equation*}
	\end{itemize}
\end{definition}




Define the Hamiltonian $H: \T^{d} \times \mathbb{R}^{d} \to \mathbb{R}$ associated with $L$ by 
$$H(x,p)=\sup_{v \in \mathbb{R}^{d}} \Big\{ \big\langle p,v \big\rangle -L(x,v) \Big\}, \quad  \forall (x,p) \in \T^{d} \times \mathbb{R}^{d}.$$
It is straightforward to check that if $L$ is a Tonelli Lagrangian, then $H$ defined above also satisfies ($i$), ($ii$), and ($iii$) in \ref{def:def1}. Such a function $H$ is called a Tonelli Hamiltonian. 

Let us recall definitions of weak KAM solutions and viscosity solutions of the Hamilton-Jacobi equation 
\begin{align}\label{eq:hj}
H(x,Du)=c, \quad x\in \T^{d},	
\end{align}
where $c$ is a real constant.


\begin{definition}[Viscosity solutions]\label{visco}
	Let $V\subset \T^{d}$ be an open set.
	\begin{itemize}
		\item [($i$)] A function $u:V\rightarrow \mathbb{R}$ is called a viscosity subsolution of equation \eqref{eq:hj}, if for every $C^1$ function $\varphi:V\rightarrow\mathbb{R}$ and every point $x_0\in V$ such that $u-\varphi$ has a local maximum at $x_0$, we have that 
		\[
		H(x_0,D\varphi(x_0))\leq c;
		\]
		\item [($ii$)] A function $u:V\rightarrow \mathbb{R}$ is called a viscosity supersolution of equation \eqref{eq:hj}, if for every $C^1$ function $\psi:V\rightarrow\mathbb{R}$ and every point $y_0\in V$ such that $u-\psi$ has a local minimum at $y_0$, we have that
		\[
		H(y_0,D\psi(y_0))\geq c;
		\]
		\item [($iii$)] A function $u:V\rightarrow\mathbb{R}$ is called a viscosity solution of equation \eqref{eq:hj} if it is both a viscosity subsolution and a viscosity supersolution.
	\end{itemize}
\end{definition}

\begin{definition}[Ma\~n\'e critical value]\label{def:def4}
The Ma\~n\'e critical value of  a Tonelli Hamiltonian $H$ is defined by
\begin{equation*} 
c(H):=\inf \left\{c \in \mathbb{R} :\ \exists\ u \in C(\T^{d}) \ \text{viscosity sol. of} \ H(x,Du) =c\right\}.
\end{equation*}
\end{definition}
See \cite[Theorem 1.1]{bib:FM} for the following weak KAM theorem for noncompact state spaces.
\begin{theorem}[Weak KAM theorem]\label{thm:wkt}
Let $H$ be a Tonelli Hamiltonian. Then, there exists a global viscosity solution of equation
\[
H(x, Du)=c(H),\quad x\in \T^{d}.
\]
\end{theorem}

We conclude this section by recalling the notion of Mather set and the role such a set plays for the regularity of viscosity solutions. Let $L$ be a Tonelli Lagrangian. As is well known, the associated Euler-Lagrange equation, i.e., 
\begin{equation}\label{eq:EL} 
\frac{d}{dt}D_{v}L(x, \dot x)=D_{x}L(x, \dot x),  
\end{equation} 
generates a flow of diffeomorphisms $\phi_{t}^{L}: \T^{d} \times \mathbb{R}^{d} \to \T^{d} \times \mathbb{R}^{d}$, with $t \in \mathbb{R}$, defined by 
\begin{equation*}\label{eq:lab11}
 \phi_{t}^{L}(x_{0},v_{0})=( x(t), \dot x(t)),\end{equation*} 
where $x: \mathbb{R} \to \T^{d}$ is the maximal solution of \eqref{eq:EL} with initial conditions $x(0)=x_{0}, \ \dot x(0)=v_{0}$. It should be noted that, for any Tonelli Lagrangian, the flow $\phi_{t}^{L}$ is complete, see for instance \cite{bib:FM}. 

 We recall that a Borel probability measure $\mu$ on $\T^{d} \times \mathbb{R}^{d}$ is called $\phi_{t}^{L}$-invariant, if $$\mu(B)=\mu(\phi_{t}^{L}(B)), \quad  \forall t \in \mathbb{R}, \quad  \forall B \in \mathscr{B}(\T^{d} \times \mathbb{R}^{d}),$$ or, equivalently,  
 \begin{equation*}
 \int_{\T^{d} \times \mathbb{R}^{d}} {f(\phi_{t}^{L}(x,v))\ \mu(dx,dv)} =\ \int_{\T^{d} \times \mathbb{R}^{d}}{f(x,v)\ \mu(dx,dv)}, \,\, \forall f \in C^{\infty}_{c}(\T^{d} \times \mathbb{R}^{d}).
 \end{equation*}
 We denote by $\mathcal{M}_{L}$ the class of all $\phi_{t}^{L}$-invariant probability measures. 
\begin{definition}[Mather measures \cite{bib:MAT}]\label{mat}
	A probability measure $\mu \in \mathcal{M}_{L}$ is called a Mather measure for $L$, if it satisfies 
	\[  \int_{\mathbb{R}^{d} \times \mathbb{R}^{d}}{L(x,v)\ \mu(dx,dv)}=\inf_{\nu \in \mathcal{M}_{L}} \int_{\mathbb{R}^{d} \times \mathbb{R}^{d}}{L(x,v)\ \nu(dx,dv)}. \]
	\end{definition} 
\noindent In \cite{bib:FA}, it was proved that 
\begin{equation*}\label{eq:lab44} c(H)=- \inf_{\nu \in \mathcal{M}_{L}} \int_{\mathbb{R}^{d} \times \mathbb{R}^{d}}{L(x,v)\ \nu(dx,dv)}.	
\end{equation*}
Denote by $\mathcal{M}_{L}^{\ast}$ the set of all Mather measures. Observe that, if $L$ (resp. $H$) is a reversible Lagrangian (resp. reversible Hamiltonian), then
\begin{equation*}\label{eq:rev}
-c(H)=\inf_{x \in \T^{d}} L(x,0).
\end{equation*}
 The Mather set is the subset $\widetilde\M_{0}\subset \T^{d}\times \mathbb{R}^d$   defined by 
\[
\widetilde\M_{0} = \overline{\bigcup_{\mu \in \mathcal{M}^{\ast}_{L}} \supp(\mu)}.
\]
We call  $\M_{0}=\pi_{1}(\mathcal{M}_0) \subset \T^{d}$ the projected Mather set. 
See \cite[Theorem 4.12.3]{bib:FA} for the following result.
\begin{theorem}\label{regularity}
	If $u$ is dominated by $L+c(H)$, then it is differentiable at every point of the projected Mather set $\M_{0}$. Moreover, if $(x,v) \in \mathcal{M}_{0}$,  then $$ Du(x) =D_{v}L(x,v) $$ and the map $M_{0} \to \T^{d} \times \mathbb{R}^{d}$, defined by $x \mapsto (x,Du(x))$, is locally Lipschitz with a Lipschitz constant which is independent of $u$.
\end{theorem}

\section{N-player games in pure strategies}
\label{sec:Purestrategies}

Fix the number of players to $N \in \N$.  For $i=1, \dots, N$ we consider a family of Lagrangians $L^{i}: \T^{d} \times \R^{d} \to \R$ satisfying the following. 
\begin{itemize}
\item[{\bf (L1)}] $L^{i}$ belongs to $C^{2}(\T^{d} \times \R^{d})$ for any $i=1, \dots, N$.
\item[{\bf (L2)}] There exists a constant $C_{0} \geq 1$ such that for any $(x_{i}, v_{i}) \in \T^{d} \times \R^{d}$
\begin{align*}
D_{v}^{2} L^{i}(x_{i}, v_{i}) \geq\ & \frac{1}{C_{0}}{\rm Id},
\\
\frac{1}{C_{0}}(1+|v_{i}|^{2}) \leq\ & L^{i}(x_{i}, v_{i}) \leq\ C_{0}(1+|v_{i}|^{2}),
\end{align*}
where ${\rm Id}$ denotes the identity matrix of dimension $d$.
\end{itemize}
Moreover, we consider a family of coupling functions 
\[
F^{i}: \T^{dN} \to \R
\]
continuous in all variables. 
Then, for any $(x_{1}, \dots, x_{N}) \in \T^{dN}$, and any $u_{i} \in L^{2}(0,T;\R^{d})$ we define the ergodic payoff
\begin{equation}\label{eq:purepayoff}
 J^{i}\big((x_{j})_{j=1}^{N}, (u_{j})_{j=1}^{N}\big)
=  \liminf_{T \to \infty} \frac{1}{T}\int_{0}^{T}\big(L^{i}(\gamma_{i}(t), u_{i}(t)) + F^{i}(\gamma_{1}(t), \dots, \gamma_{N}(t)) \big)\ {\rm dt}
\end{equation}
where $\gamma_{i}$ satisfies the state equation 
\begin{equation*}
\dot\gamma_{i}(t)=u_{i}(t), \quad \gamma_{i}(0)=x_{i}, 
\end{equation*}
for any $i \in \{1, \dots, N\}$. 

Next, we give the definition of Nash equilibria in pure strategies. 
\begin{definition}[{\bf Nash equilibria in pure strategies}]
	We say that $(\overline{u}_{j})_{j=1}^{N} \in L^{2}(0,T; \R^{d})$ is a {\em Nash equilibrium in pure strategies} for the initial condition $(x_{j})_{j=1}^{N} \in \T^{dN}$ if 
	\begin{equation*}
	J^{i}\big(	(x_{j})_{j=1}^{N},(\overline{u}_{j})_{j=1}^{N} \big) \leq J^{i}\big((x_{j})_{j=1}^{N}, (\overline{u}_{j})_{j\not=i}, u_{i}\big)
	\end{equation*}
for any control $u_{i} \in L^{2}(0,T;\R^{d})$. 
\end{definition}

In this first part we are interested in finding Nash equilibria in pure strategies for \eqref{eq:purepayoff}. To do so, we consider the system of Hamilton-Jacobi equations associated with the game \eqref{eq:purepayoff}, that is,
\begin{equation}\label{eq:Nhj}
	H^{i}(x_{i}, D_{x_{i}}v^{i}) + \sum_{j\not= i}{\langle D_{p}H^{j}(x_{j}, D_{x_{j}}v^{j}), D_{x_{j}}v^{j}\rangle} 
	+ \lambda_{i} = F^{i}(x_{1}, \dots, x_{N}), \,\, \text{on}\,\, \T^{dN}
\end{equation}
where $H^{k}$ is defined by
\begin{equation*}
H^{k}(x,p)= \sup_{u \in \R^{d}}\big\{\langle p, u\rangle - L^{k}(x, u)\big\} \quad (k \in \{1, \dots, N\}).
\end{equation*}
We will show that the existence of Nash equilibria in pure strategies for \eqref{eq:purepayoff} depends on the existence of smooth solutions to \eqref{eq:Nhj}. In view of \eqref{eq:purepayoff} and \eqref{eq:Nhj}, inspired by weak KAM theory for Tonelli Hamiltonian systems, we expect that if $(\overline{u}_{1},\dots, \overline{u}_{N})$ is a Nash equilibrium then 
\begin{equation*}
J^{i}((x_{j})_{j=1}^{N}, (\overline{u}_{j})_{j=1}^{N})=\lambda_{i} \,\, (i=1,\dots, N). 
\end{equation*}

\medskip
{\it Hereafter, we assume that there exists a classical solution $\overline{v}^{i}$, for any $i=1,\dots,N$ to \eqref{eq:Nhj}. Then, having such a solution at our disposal we construct a Nash equilibrium in pure strategies for \eqref{eq:purepayoff}. Clearly, since the existence of such a smooth solution is not always guaranteed this result will not be satisfactory. This force us to modify the set of strategies for the game, i.e., we look for strategies on the space of Wasserstein probability measures.}
\medskip

\begin{remarks}\em
It is well-know that, for instance, there exists a smooth solution to Hamilton-Jacobi equations if the Lagrangian function of the underlying minimization problem is convex not only with respect to the velocity variable but also with respect to the state variable, see for instance \cite{Rockafellar_2000, Rockafellar_2001}. That is, besides the assumptions {\bf (L1)} and {\bf (L2)} if we require $L$ to satisfy
\begin{equation*}
D_{x}^{2} L^{i}(x_{i}, v_{i}) \geq\  \frac{1}{C}{\rm Id},
\end{equation*}
for some constant $C > 0$ we obtain that \eqref{eq:Nhj} has a smooth solution. Moreover, we stress that the existence of such smooth solutions is not a necessary condition for the existence of Nash equilibria in pure strategies but it is sufficient to obtain also a representation of that. 
\end{remarks}

Next, through a series of Lemmas we proceed with the construction of Nash equilibria in pure strategies.


\begin{lemma}
Fix $i \in \{1, \dots, N\}$ and let $(\overline{v}^{i}, \overline\lambda_{i})$ be a solution to \eqref{eq:Nhj} such that $\overline{v}^{i}$ is smooth. Let $a$, $b \in \R$ with $a < b$. Then, for any family of absolutely continuous curves 
\[ 
\gamma_{j}:[a,b] \to \T^{d} \quad  (j=1,\dots,N)
\] 
we have that 
\begin{multline*}
 \overline{v}^{i}(\gamma_{1}(b), \dots, \gamma_{N}(b)) - \overline{v}^{i}(\gamma_{1}(a), \dots, \gamma_{N}(a)) 
	\leq\ \int_{a}^{b}{\big(L^{i}(\gamma_{i}(t), u_{i}(t)) + F^{i}(\gamma_{1}(t), \dots, \gamma_{N}(t)) \big)\ {\rm dt}}  
	\\
	+ \overline\lambda_{i}(b-a).
\end{multline*}
\end{lemma}
\proof
Given $a$, $b \in \R$ with $a < b$, fix $j \in \{1, \dots, N\}$ and take $\gamma_{j}:[a,b] \to \T^{d}$ with $\gamma_{j} \in \text{AC}([0,T]; \T^{d})$. Then, we have that 
\begin{align*}
 \overline{v}^{i}(\gamma_{1}(b), \dots, \gamma_{N}(b)) - \overline{v}^{i}(\gamma_{1}(a), \dots, \gamma_{N}(a))
=\ 	\int_{a}^{b}{\sum_{j=1}^{d}{\langle D_{x_{j}}v^{j}(\gamma_{1}(t), \dots, \gamma_{N}(t)), \dot\gamma_{j}(t))\rangle}\  {\rm dt}}.
\end{align*}
Hence, by definition of Legendre Transform applied to the $i$-th term of the sum on the right hand side we obtain
\begin{multline*}
	 \overline{v}^{i}(\gamma_{1}(b), \dots, \gamma_{N}(b)) - \overline{v}^{i}(\gamma_{1}(a), \dots, \gamma_{N}(a))
=\ \\ \int_{a}^{b}{\sum_{j=1}^{d}{\langle D_{x_{j}}v^{j}(\gamma_{1}(t), \dots, \gamma_{N}(t)), \dot\gamma_{j}(t))\rangle}\  {\rm dt}}
\\
\leq\  \int_{a}^{b}\Big(\sum_{j\not= i}{\langle D_{x_{j}}v^{j}(\gamma_{1}(t), \dots, \gamma_{N}(t)), \dot\gamma_{j}(t))\rangle} 
\\
+\  L^{i}(\gamma_{i}(t), u_{i}(t)) + H^{i}(\gamma_{i}(t), D_{x_{i}}v^{i}(\gamma_{1}(t), \dots, \gamma_{N}(t)))\Big)\  {\rm dt}.
\end{multline*}
Finally, since $\overline{v}^{i}$ solves \eqref{eq:Nhj} we get 
\begin{multline*}
	\overline{v}^{i}(\gamma_{1}(b), \dots, \gamma_{N}(b)) - \overline{v}^{i}(\gamma_{1}(a), \dots, \gamma_{N}(a))
	\leq\  \int_{a}^{b}{\big(L^{i}(\gamma_{i}(t), u_{i}(t)) + F^{i}(\gamma_{1}(t), \dots, \gamma_{N}(t)) \big)\ {\rm dt}} \\ + \overline\lambda_{i}(b-a)
\end{multline*}
and this yields the conclusion. \qed

This suggest the following definition of dominated functions and that of calibrated curves.

\begin{definition}[{\bf Dominated functions}]
	Let $\varphi: \T^{dN} \to \R$ be a continuous function and let $a$, $b$, $c \in \R$ be such that $a < b$. Then, $\varphi$ is said to be dominated by $L^{i} + F^{i} + c$  (with $i=1,\dots,N$), which we denote by $\varphi \prec L^{i} + F^{i} + c$, if for any family of absolutely continuous curves $(\gamma_{j})_{j=1}^{N}$ with $\gamma_{j}:[a,b] \to \T^{d}$ we have that 
	\begin{align}\label{eq:dominated}
	\begin{split}
	& \varphi(\gamma_{1}(b), \dots, \gamma_{N}(b)) - \varphi(\gamma_{1}(a), \dots, \gamma_{N}(a))
	\\
	\leq\ & \int_{a}^{b}{\big(L^{i}(\gamma_{i}(t), u_{i}(t)) + F^{i}(\gamma_{1}(t), \dots, \gamma_{N}(t)) \big)\ {\rm dt}} + c(b-a).
	\end{split}
\end{align}
\end{definition}

\begin{definition}[{\bf Calibrated curves}]
Let $\varphi: \T^{dN} \to \R$ be a continuous function and let $a$, $b$, $c \in \R$ be such that $a < b$. Then, the family of absolutely continuous curves $(\gamma_{j})_{j=1}^{N}$ with $\gamma_{j}:[a,b] \to \T^{d}$ is said to be calibrated for $(\varphi, c)$ if  
\begin{multline*}
	 \varphi(\gamma_{1}(b), \dots, \gamma_{N}(b)) - \varphi(\gamma_{1}(a), \dots, \gamma_{N}(a))
	=\  \int_{a}^{b}{\big(L^{i}(\gamma_{i}(t), u_{i}(t)) + F^{i}(\gamma_{1}(t), \dots, \gamma_{N}(t)) \big)\ {\rm dt}} \\ + c(b-a).
\end{multline*}
\end{definition}

\begin{theorem}\label{mainpure}
	Fix $i \in \{1, \dots, N\}$ and let $(\overline{v}^{i}, \overline\lambda_{i})$ be a solution to \eqref{eq:Nhj} such that $\overline{v}^{i}$ is smooth. Then, for any $i \in \{1,\dots, N\}$ the control
\begin{equation*}
\overline{u}_{i}(x)=D_{p}H^{i}(x, D\overline{v}^{i}(x))
\end{equation*}
is a Nash equilibrium in pure strategies and
\[
J^{i}\big( (x_{j})_{j=1}^{N}, (\overline{u})_{j=1}^{N} \big)=\overline\lambda_{i}
\] 
for any initial condition $(x_{j})_{j=1}^{N} \in \T^{dN}$. 
\end{theorem}
\proof 
We first show that the trajectory associated with $\overline{u}_{i}(x)=D_{p}H^{i}(x, D\overline{v}^{i}(x))$, denote by $\overline\gamma_{i}$ is calibrated for $(\overline{v}^{i}, \overline\lambda_{i})$ for any $i=1,\dots, N$. Indeed, since $\overline{v}^{i}$ is smooth and $\overline\gamma_{i}$ is absolutely continuous from the definition of Legendre Transform we have that 
\begin{align*}
& \frac{1}{T}\big(\overline{v}^{i}(\overline\gamma_{1}(T), \dots, \overline\gamma_{N}(T))- 	\overline{v}^{i}(x_{1}, \dots, x_{N})\big) 
\\
=\ & \frac{1}{T}\int_{0}^{T}{\big(L^{i}(\gamma_{i}(s), u_{i}(s)) + F^{i}(\overline\gamma_{1}(s), \dots, \overline\gamma_{N}(s)) \big)\ {\rm ds}} - \overline\lambda_{i}. 
\end{align*}
Hence, by compactness of $\T^{d}$, as $T \to \infty$ we get 
\begin{equation*}
J^{i}\big( (x_{j})_{j=1}^{N}, (\overline{u}_{j})_{j=1}^{N}\big)=\overline\lambda_{i}
\end{equation*}
for any $i=1,\dots,N$.

Next, given any other strategy $u_{i} \in L^{2}(0,T; \R^{d})$ with $u_{i}=\overline{u}_{i}$ (for $i=1,\dots,N$) the curve $\gamma_{i}$, generated by $u_{i}$, satisfies  \eqref{eq:dominated} for $\varphi=\overline{v}^{i}$ and $c=\overline\lambda_{i}$. Hence, we obtain 
\begin{multline*}
 \frac{1}{T}\big(v^{i}(\overline\gamma_{1}(T), \dots, \overline\gamma_{i-1}(T), \gamma_{i}(T), \overline\gamma_{i+1}(T), \dots, \overline\gamma_{N}(T)) - 	v^{i}(x_{1}, \dots, x_{N})\big) 
\leq\  \frac{1}{T}\int_{0}^{T} \big(L^{i}(\gamma_{i}(s), u_{i}(s)) \\ + F^{i}(\overline\gamma_{1}(s), \dots, \overline\gamma_{i-1}, \gamma_{i}(s), \overline\gamma_{i+1}(s), \overline\gamma_{N}(s)) \big)\ {\rm ds} - \overline\lambda_{i}
\end{multline*}
and as $T \to \infty$ we get
\begin{equation*}
\overline\lambda_{i} \leq J^{i}\big( (x_{j})_{j=1}^{N}, (\overline{u}_{j})_{j\not=i}, u_{i}). \eqno\square	
\end{equation*}


\section{From Mather theory to mixed strategies}
\label{sec:Mathertheory}

In this section we exploit the connections between the N-players game and the Mather theory for Tonelli Hamiltonian systems. This allow us, in the next section, to define mixed strategies and establish the connection between Nash equilibria in mixed strategies and Mather measures. 

 We want to derive from \eqref{eq:Nhj} a different system of PDEs from which we can synthesize a Nash equilibrium and from which we can find an intuitive connection between Mather theory and Nash equilibria. To do so,

\medskip
\qquad\qquad\qquad{\it Assume that $D_{x_{i}} v^{i}$ depends only on $x_{i}$ for any $i=1,\dots, N$.}
\medskip

Then, equations in \eqref{eq:Nhj} are decoupled and we can derive a PDE for the candidate measure that will define mixed strategies for the game. 
A game that satisfies the above assumption is the one with independent players and payoff, for instance, 
\begin{equation*}
\liminf_{T \to \infty} \frac{1}{T}\int_{0}^{T}\big(L^{i}(\gamma_{i}(t), u_{i}(t)) + F^{i}(\gamma_{i}(t)) \big)\ {\rm dt}.
\end{equation*}
However, observe that the above assumption is used only in this section to give an intuitive roadmap toward the use of Mather theory and its connection with the differential N-player game.

\begin{proposition}\label{prop:continuityeq}
	Fix $i \in \{1, \dots, N\}$ and let $(\overline{v}^{i}, \overline\lambda_{i})$ be a solution to \eqref{eq:Nhj} such that $\overline{v}^{i}$ is smooth. Let $\mu^{i}$ be an invariant measure for the Euler flow $\Phi_{t}^{i}(\cdot, \cdot)$ of the Lagrangian $L^{i}(x_{i}, u_{i}) + F^{i}(x_{1}, \dots, x_{N})$. Then, the probability measure $m^{i}= \pi_{1} \sharp \mu^{i}$ solves
	\begin{equation}\label{eq:stationaryeq}
	\ddiv\big(m^{i}D_{p}H^{i}(x_{i}, D_{x_{i}}\overline{v}^{i}) \big)=0,
	\end{equation}
	in the sense of distributions.
\end{proposition}
\proof
Let $(\gamma^{i}, \dot\gamma^{i})$ denote the Euler flow associated with the Lagrangian $L^{i}(x_{i}, u_{i}) + F^{i}(x_{1}, \dots, x_{N})$, i.e., $\Phi_{t}^{i}(x, v)=(\gamma^{i}, \dot\gamma^{i})$ and
\begin{equation*}
\begin{cases}
{\rm \frac{d}{dt}} D_{v} L^{i}(\gamma^{i}(t), \dot\gamma^{i}(t)) = D_{x}\big( L^{i}(\gamma^{i}(t), \dot\gamma^{i}(t)) \\ +\,\,  F^{i}(\gamma^{1}(t), \dots, \gamma^{N}(t)) \big), \,\, (t \geq 0)
\\
\gamma^{i}(0)=x, \quad \dot\gamma^{i}(0)=v.
\end{cases}
\end{equation*}
Moreover, since $\mu^{i}$ is invariant for $(\gamma^{i}, \dot\gamma^{i})$ then the probability measure $m^{i}= \pi_{1} \sharp \mu^{i}$ is invariant for $\gamma^{i}$. Hence, for any $\psi \in C^{1}(\T^{d})$ we get
\begin{multline*}
0 = {\rm \frac{d}{dt}} \int_{\T^{d}}{\psi(\gamma^{i}(t))\ m^{i}({\rm dx_{i}})}   =\  \int_{\T^{d}}{\langle D\psi(\gamma^{i}(t)), D_{p}H^{i}(\gamma^{i}(t), D_{x_{i}}\overline{v}^{i}(\gamma^{i}(t)) \rangle\ m^{i}({\rm dx_{i}})}
\\
=\  \int_{\T^{d}}{\langle D\psi(x_{i}), D_{p}H^{i}(x_{i}, D_{x_{i}}\overline{v}^{i}(x_{i})) \rangle\ m^{i}({\rm dx_{i}})}
\end{multline*}
which completes the proof. \qed

\medskip

\smallskip
Define the function $z^{i}: \T^{d} \to \R$ as
\begin{equation}\label{eq:singleplayer}
z^{i}(x)=v^{i}(x_{1}, \dots, x_{i-1}, x,x_{i+1}, x_{N}), \quad (x \in \T^{d})	
\end{equation}
and, given $(m_{j})_{j=1}^{N} \in \PP(\T^{d})$ define the function $V^{i}[m_{1}, \dots, m_{N}](\cdot): \T^{d} \to \R$ by
\begin{equation}\label{eq:VI}
V^{i}[m_{1}, \dots, m_{N}](y) = \int_{T^{d(N-1)}}{F^{i}(x_{1}, \dots, x_{i-1}, y, x_{i+1}, \dots, x_{N})\ \prod_{j\not=i}{ m^{j}({\rm dx_{j})}}}.	
\end{equation}

 \begin{proposition}
 	Fix $i \in \{1, \dots, N\}$, let $(\overline{v}^{i}, \overline\lambda_{i})$ be a smooth solution of \eqref{eq:Nhj} and let $z^{i}$ be defined as in \eqref{eq:singleplayer}. Let $(\mu^{j})_{j=1}^{N}$ be an invariant probability measure for the Euler flow of $(\Phi_{t}^{j})_{j=1}^{N}$. Then, $z^{i}$ is a viscosity solution of 
 	\begin{equation}\label{eq:1HJ}
 		H^{i}(x, Dz^{i}(x))+\overline\lambda_{i}=V^{i}[m_{1}, \dots, m_{N}](x) \quad (x \in \T^{d})
 	\end{equation}
 	where $(m_{j})_{j=1}^{N}=(\pi_{1} \sharp \mu^{j})_{j=1}^{N}$. 
 \end{proposition}
\proof
We fix $i \in \{1, \dots, N\}$, multiply the equation \eqref{eq:Nhj} by $\prod_{j \not= i} m_{j}({\rm dx_{j})}$ and integrate over $\T^{d(N-1)}$. Observe that, for $k \not= i$ from \eqref{eq:stationaryeq} we get
\begin{align*}
	 & \int_{\T^{d(N-1)}}{\langle D_{p}H^{i}(x_{k}, D_{x_{k}}\overline{v}^{k}), D_{x_{k}}\overline{v}^{k}\rangle \ \prod_{j \not= i}m_{j}({\rm dx^{j})}}
	 \\
	=\ &  \int_{\T^{d(N-2)}}\int_{\T^{d}}\langle D_{p}H^{i}(x_{k}, D_{x_{k}}\overline{v}^{k}), D_{x_{k}}\overline{v}^{k}\rangle\ m_{k}({\rm dx^{k}}) \prod_{j \not= i, k}m_{j}({\rm dx^{j})} = 0.
\end{align*}
Then, since $D_{x_{i}}\overline{v}^{i}$ depends only on $x_{i}$ we deduce that 
\begin{equation*}
\int_{\T^{d(N-1)}}{\big(H^{i}(x_{i}, D_{x_{i}}\overline{v}^{i}) + \overline\lambda_{i} \big)\ \prod_{j \not= i}m_{j}({\rm dx^{j}})} = H^{i}(x_{i}, D_{x_{i}}\overline{v}^{i}) + \overline\lambda_{i}.
\end{equation*}
Hence, combining together the above expressions we get \eqref{eq:1HJ}. \qed

Equation \eqref{eq:1HJ} coupled with \eqref{eq:stationaryeq} define the N-players mean field game (MFG) system 
\begin{equation}\label{eq:Nplayersystem}
\begin{cases}
	H^{i}(x, Dz^{i}(x))+\overline\lambda_{i}=V^{i}[m_{1}, \dots, m_{N}](x), & x \in \T^{d}
	\\
	\ddiv\big(m^{i}D_{p}H^{i}(x_{i}, D_{x_{i}}v^{i}) \big)=0, & x \in \T^{d}.
\end{cases}
\end{equation}

As we explained so far, the construction of Nash equilibria in pure strategies in \ref{mainpure} is not satisfactory since it depends on the existence of smooth solutions to \eqref{eq:Nhj} which, in general, is not guaranteed. Such reason, together with system \eqref{eq:Nplayersystem} motivated us to introduce the notion of mixed strategies as probability measures on $\T^{d}$. In view of \ref{mainpure} we know that 
\begin{equation*}
J^{i}\big((x_{j})_{j=1}^{N}, (\overline{u}_{j})_{j=1}^{N}\big)=\overline\lambda_{i}.
\end{equation*}
where $(\overline{u}_{j})_{j=1}^{N}$ is a Nash equilibrium in pure strategies. Therefore, we expect that the value of the game is not lower than $\overline\lambda_{i}$ even playing with mixed strategies. 

\section{N-player game in mixed strategies}
\label{sec:Mixedstrategies}

\subsection{Heuristic motivation}
{\it Here, we give the heuristic motivation that inspired us to give the definition of the new payoff with strategies given by probability measures.} 

Set
\begin{equation*}
\II = \big\{(\mu^{j})_{j=1}^{N} \in \PP(\T^{d}): \mu^{j}\,\, \text{is invariant for}\,\, \Phi_{t}^{j}\big\}.
\end{equation*}
Let $(\overline{u}_{j})_{j=1}^{N}$ be a Nash equilibrium in pure strategies and consider the payoff \eqref{eq:purepayoff} evaluated at the Nash equilibrium in pure strategies given by 
\begin{align}\label{eq:heu}
\liminf_{T \to \infty} \frac{1}{T}\int_{0}^{T}{\big(L^{i}(\overline\gamma_{i}(s), \overline{u}_{i}(s)) + F^{i}(\overline\gamma_{1}(s), \dots, \overline\gamma_{N}(s)) \big)\ {\rm ds}}.
\end{align}
Let $\mu = \prod_{j=1}^{d} \mu^{j}$ with $(\mu^{j})_{j=1}^{N} \in \II$ and multiply \eqref{eq:heu} by such measure. So, we obtain 
\begin{equation*}
\liminf_{T \to \infty} \frac{1}{T}\int_{0}^{T}\int_{\T^{dN} \times \R^{dN}} \big(L^{i}(\overline\gamma_{i}(s), \dot{\overline\gamma}_{i}(s)) + F^{i}(\overline\gamma_{1}(s), \dots, \overline\gamma_{N}(s)) \big)\ \prod_{j=1}^{d} \mu^{j}({\rm dx_{j}, dv_{j}}){\rm ds}.
\end{equation*}
So, from the invariance of $(\mu^{j})_{j=1}^{N}$ w.r.t. $(\Phi^{j}_{t})_{j=1}^{N}$ we get
 \begin{align*}
	\int_{\T^{dN} \times \R^{dN}}\big(L^{i}(x_{i}, v_{i}) + F^{i}(x_{1}, \dots, x_{N}) \big)\ \prod_{j=1}^{d} \mu^{j}({\rm dx_{j}, dv_{j}})
 \end{align*}
which by definition of $V^{i}$ in \eqref{eq:VI} can be written as
\begin{equation*}
\int_{\T^{d} \times \R^{d}} \big(L^{i}(x_{i}, v_{i}) + V^{i}[\pi\sharp\mu^{1}, \dots, \pi\sharp\mu^{N}](x_{i}) \big)\ \mu^{i}({\rm dx_{i}, dv_{i}}).
\end{equation*}
This is, indeed, the functional we will consider as payoff of our game in mixed strategies.

\subsection{Existence of Nash equilibria in mixed strategies}

Now, we can proceed to the analysis of the N-player game in mixed strategies. Define the payoff function $J^{i}: \PP(\T^{d})^N \times \PP(\T^{d})^N \to \R$
\begin{equation*}
J^{i}\big((\mu_{j})_{j=1}^{N}, (\eta^{j})_{j=1}^{N}\big) = \int_{\T^{dN} \times \R^{dN}}{\big(L^{i}(x_{i}, v_{i}) + V[\pi\sharp\mu_{1}, \dots, \pi\sharp\mu_{N}](x_{i}) \big)\ \eta^{i}({\rm dx_{i}dv_{i}})}.
\end{equation*}

\begin{definition}[{\bf Nash equilibria in mixed strategies}]
	We say that $(\overline\mu^{j})_{j=1}^{N} \in \II$  is a {\em Nash equilibrium in mixed strategies} if 
\begin{equation*}
J^{i}\big((\overline\mu_{j})_{j=1}^{N}, (\overline\mu_{j})_{j=1}^{N}\big) \leq J^{i}\big((\overline\mu_{j})_{j=1}^{N}, (\overline\mu_{j})_{j\not=i}, \mu^{i}\big)
\end{equation*}
for any $\mu^{i}$ invariant for $\Phi^{i}_{t}$.
\end{definition}

Note that, equivalently, $(\overline\mu^{j})_{j=1}^{N} \in \II$ is a Nash equilibrium in mixed strategies if
\begin{align*}
& \int_{\T^{d} \times \R^{d}} \big(L^{i}(x_{i}, v_{i}) + V^{i}[\pi \sharp \overline\mu_{1}, \dots, \pi \sharp \overline\mu_{N}](x_{i}) \big)\ \overline\mu^{i}({\rm dx_{i}, dv_{i}})
\\
\leq\ &  \int_{\T^{d} \times \R^{d}} \big(L^{i}(x_{i}, v_{i}) + V^{i}[\pi \sharp \overline\mu_{1}, \dots, \pi \sharp \overline\mu_{N}](x_{i}) \big)\ \mu^{i}({\rm dx_{i}, dv_{i}})
\end{align*}
	for any $\mu^{i}$ invariant for $\Phi^{i}_{t}$.

\begin{theorem}\label{thm:mixedstrat}
There exists at least one Nash equilibrium in mixed strategies. 	
\end{theorem}
\proof
The proof is based on the application of the Kakutani fixed-point theorem. To do so, we define the set-valued map
\begin{equation*}
\E: (\PP(\T^{d})^{N}, d_{1}) \rightrightarrows  (\PP(\T^{d})^{N}, d_{1})
\end{equation*}
such that 
\begin{equation*}
\E \left((m_{j})_{j=1}^{N}\right)=\left\{(\pi \sharp \nu_{i})_{i=1}^{N}: (\nu_{i})_{i=1}^{N} \in M((m_{j})_{j=1}^{N}) \right\}
\end{equation*}
where $M((m_{j})_{j=1}^{N}) \subset \PP(\T^{d} \times \R^{d})^{N}$ is the set of minimizers of the map 
\begin{equation}\label{eq:measuremin}
\eta \mapsto \int_{\T^{d} \times \R^{d}}{\big(L^{i}(x_{i}, v_{i}) + V[m_{1}, \dots, m_{N}](x_{i}) \big)\ \eta({\rm dx_{i}dv_{i})}}
\end{equation} 
over the set $\II$. That is, given $(\nu_{i})_{i=1}^{N} \in M((m_{j})_{j=1}^{N})$ we have that 
\begin{align*}
& \inf_{ \eta \in \II} \int_{\T^{d} \times \R^{d}}{\big(L^{i}(x_{i}, v_{i}) + V[m_{1}, \dots, m_{N}](x_{i}) \big)\ \eta({\rm dx_{i}dv_{i})}}
\\
= & \int_{\T^{d} \times \R^{d}}{\big(L^{i}(x_{i}, v_{i}) + V[m_{1}, \dots, m_{N}](x_{i}) \big)\ \nu_{i}({\rm dx_{i}dv_{I})}}.
\end{align*}
For any $(m_{j})_{j=1}^{N} \in \PP(\T^{d})$ we have that $\E((m_{j})_{j=1}^{N})$ is non-empty, convex and compact set. Hence, we have to check that $\E$ has closed graph in order to apply the Kakutani fixed-point Theorem. 

Fix a sequence $(m^{k}_{j})_{j=1}^{N} \in \PP(\T^{d})$ and $(\nu^{k}_{i})_{i=1}^{N} \in \PP(\T^{d})$ for $k \in \N$ such that 
\begin{equation*}
m^{k}_{j} \to \overline{m}_{j}, \quad \nu^{k}_{i} \to \overline\nu_{i} \quad \text{as}\,\, k \to \infty\,\, \text{w.r.t.}\,\, d_{1}, \text{and} \quad (\nu^{k}_{i})_{i=1}^{N} \in \E((m^{k}_{j})_{j=1}^{N}) \,\, \forall\ j \in \N. 
\end{equation*}
For any $k \in \N$, let $(\eta^{k}_{i})_{i=1}^{N} \in \II$ be such that $(\nu_{i}^{k})_{i=1}^{N}=(\pi \sharp \eta^{k}_{i})_{i=1}^{N}$. For simplicity of notation, fix $i \in \{1, \dots, N\}$. From the invariance of $\nu_{i}^{k}$ w.r.t. the Euler flow and the compactness of the latter \cite[Corollary 4.4.5]{bib:FA} we have that $\{\eta_{i}^{k}\}_{k \in \N}$ is tight and by Prokhorov's Theorem there exists $\overline\eta_{i}$ such that $\eta_{i}^{k} \to \overline\eta_{i}$ as $k \to \infty$ and $\overline\nu_{i}= \pi \sharp \overline\eta_{i}$. We have to show that $\overline\eta_{i}$ is a minimizing measure associated with $(\overline{m}_{j})_{j=1}^{N}$. 

For simplicity of notation, we set
\begin{equation*}
\Lambda((m_{j})_{j=1}^{N})  = \min_{\nu \in \K} \int_{\T^{d} \times \R^{d}}{\big(L^{i}(x_{i}, v_{i}) + V^{i}[m_{1}, \dots, m_{N}](x_{i}) \big)\ \nu({\rm dx_{i}dv_{i})}}.
\end{equation*}
From the lower-semicontinuity of the function $V$ w.r.t. the measure variable, we have that 
\begin{multline}\label{eq:lowersemi}
 \int_{\T^{d} \times \R^{d}}{\big(L^{i}(x_{i}, v_{i}) + V^{i}[\overline{m}_{1}, \dots \overline{m}_{N}](x_{i}) \big)\ \overline\eta_{i}({\rm dx_{i}dv_{i}})}
\\
\leq\  \liminf_{k \to \infty} \int_{\T^{d} \times \R^{d}}\big(L^{i}(x_{i}, v_{i}) + V^{i}[m^{k}_{1}, \dots, m^{k}_{N}](x_{i}) \big)\ \eta^{k}_{i}({\rm dx_{i}dv_{i}}). 
\end{multline}
Thus, we proceed to prove that the right-hand side of \eqref{eq:lowersemi} is not larger than $\Lambda((\overline{m}_{j})_{j=1}^{N})$.  Let $(\widetilde\nu_{i})_{i=1}^{N} \in \E((\overline{m}_{j})_{j=1}^{N})$, let $(\widetilde\eta_{i})_{i=1}^{N}$ be such that $(\widetilde\nu_{i})_{i=1}^{N}=(\pi\sharp\widetilde\eta_{i})_{i=1}^{N}$ and fix $\eps > 0$. Then, there exists $R \geq 0$ such that 
\begin{equation*}
\int_{\T^{d} \times \R^{d}\backslash B_{R}}{|v_{i}|^{2}\ \widetilde\eta_{i}({\rm dx_{i}dv_{i}})} \leq \eps, \quad \forall\ i \in \{1, \dots, N\}. 
\end{equation*}
Fix $i \in \{1, \dots, N\}$. Then, since $(\nu^{k}_{i})_{i=1}^{N}=(\pi \sharp \eta^{k}_{i})_{i=1}^{N} \in \E((m^{k}_{j})_{j=1}^{N})$ we get
\begin{multline*}
\Lambda((m^{k}_{j})_{j=1}^{N})  =\ \int_{\T^{d} \times \R^{d}}{\big(L^{i}(x_{i}, v_{i}) + V^{i}[m^{k}_{1}, \dots, m^{k}_{N}](x_{i}) \big)\ \eta^{k}_{i}({\rm dx_{i}dv_{I}})}
\\
\leq\  \int_{\T^{d} \times \R^{d}}{\big(L^{i}(x_{i}, v_{i}) + V^{i}[m^{k}_{1}, \dots, m^{k}_{N}](x_{i}) \big)\ \widetilde\eta_{i}({\rm dx_{i}dv_{i}})}
\\
\leq\  \int_{\T^{d} \times \overline{B}_{R}}{\big(L^{i}(x_{i}, v_{i}) + V^{i}[m^{k}_{1}, \dots, m^{k}_{N}](x_{i}) \big)\ \widetilde\eta_{i}({\rm dx_{i}dv_{i}})} \\
+\  \int_{\T^{d} \times \R^{d}\backslash B_{R}}{\big(L^{i}(x_{i}, v_{i}) + V^{i}[m^{k}_{1}, \dots, m^{k}_{N}](x_{i}) \big)\ \widetilde\eta_{i}({\rm dx_{i}dv_{i}})}
\\
\leq\  \int_{\T^{d} \times \overline{B}_{R}}{\big(L^{i}(x_{i}, v_{i}) + V^{i}[m^{k}_{1}, \dots, m^{k}_{N}](x_{i}) \big)\ \widetilde\eta_{i}({\rm dx_{i}dv_{i}})} 
\\
+\  \int_{\T^{d} \times \R^{d}\backslash B_{R}}\Big(C_{0}(1+|v_{i}|^{2}) + \|V^{i}[m^{k}_{1}, \dots, m^{k}_{N}](\cdot)\|_{\infty}\Big)\ \widetilde\eta_{i}({\rm dx_{i}dv_{i}})
\\
\leq\  \int_{\T^{d} \times \overline{B}_{R}}\big(L^{i}(x_{i}, v_{i}) + V^{i}[m^{k}_{1}, \dots, m^{k}_{N}](x_{i}) \big)\ \widetilde\eta_{i}({\rm dx_{i}dv_{i}}) \\
+\  \int_{\T^{d} \times \R^{d}\backslash B_{R}}\Big(C_{0}(1+|v_{i}|^{2}) + \|V^{i}[\overline{m}_{1}, \dots, \overline{m}_{N}](\cdot)\|_{\infty}\Big)\ \widetilde\eta_{i}({\rm dx_{i}dv_{i}})
\\
\leq\  \Lambda(\mu) + 2\eps. 
\end{multline*}
Combining this inequality with \eqref{eq:lowersemi}, by arbitrarily of $\eps$ we obtain 
\begin{equation*}
\int_{\T^{d} \times \R^{d}}{\big(L^{i}(x_{i}, v_{i}) + V^{i}[\overline{m}_{1}, \dots \overline{m}_{N}](x_{i}) \big)\ \overline\eta_{i}({\rm dx_{i}dv_{i}})} \leq \Lambda(\mu).
\end{equation*}
Iterating the same argument for any $i \in \{1, \dots, N\}$, we get $(\overline\nu_{i})_{i=1}^{N} \in \E((\overline{m})_{j=1}^{N})$. Therefore, by Kakutani fixed-point Theorem there exists $(\overline{p}_{j})_{j=1}^{N} \in \E \left((\overline{p}_{j})_{j=1}^{N}\right)$, which implies that $(\overline\mu_{j})_{j=1}^{N} \in M((\overline{p}_{j})_{j=1}^{N})$ such that $(\overline{p}_{j})_{j=1}^{N}=(\pi \sharp \overline\mu_{j})_{j=1}^{N}$ is a Nash equilibrium in mixed strategies. \qed


\begin{corollary}
Let $(\overline\mu^{j})_{j=1}^{N} \in \II$  be a Nash equilibrium in mixed strategies. Then, we have that 
\begin{equation*}
	J^{i}((\overline\mu^{j})_{j=1}^{N}, (\overline\mu^{j})_{j=1}^{N})=\lambda_{i}.
\end{equation*}
\end{corollary}

\section{Limit of Nash equilibria}
\label{sec:Meanfieldlimit}

In this section, we consider the case of symmetric game and we study the limit of the N-player games as $N$ goes to infinity. To do so, we assume that players are indistinguishable and, consequently, on the data of the model we assume that 
\begin{itemize}
\item[({\bf S})] $L^{i}=L^{j}$, and $F^{i}=F^{j}$ for any $i, j \in \{1, \dots, N\}$. 
\end{itemize}
In particular, this allowed us to derive the ergodic PDEs system used in \cite{bib:CAR}, \cite{bib:CCMW, bib:CCMW1} to describe the long-time behavior of solutions to the classical MFG system with finite horizon $T$. We recall that, a function $\phi_{N}: \T^{dN} \to \R$ is a symmetric function if
\begin{equation*}
\phi_{N}(x_{1}, \dots, x_{N})=\phi(x_{\sigma(1)}, \dots, x_{\sigma(N)}) 
\end{equation*}
for all permutation $\sigma$ on $\{1, \dots, N\}$. 

Next, from the same reasoning to those in Theorem \ref{thm:mixedstrat} we get the following result.

\begin{proposition}
Assume {\bf (S)}. Then, there exists a Nash equilibrium in mixed strategies $(\overline\mu_{j})_{j=1}^{N}$ such that $\overline\mu_{j}=\overline\mu_{k}$ for any $j, k \in \{1, \dots, N\}$. 
\end{proposition}

\proof 
Define the set-valued map
\begin{equation*}
\E: (\PP(\T^{d}), d_{1}) \rightrightarrows  (\PP(\T^{d}), d_{1})
\end{equation*}
such that 
\begin{equation*}
\E (m)=\left\{\pi \sharp \nu: \nu \in M(m) \right\}
\end{equation*}
where $M(m) \subset \PP(\T^{d} \times \R^{d})$ is the set of minimizers of the map 
\begin{equation}\label{eq:measuremin}
\eta \mapsto \int_{\T^{d} \times \R^{d}}{\big(L^{i}(x_{i}, v_{i}) + V[m, \dots, m](x_{i}) \big)\ \eta({\rm dx_{i}dv_{i})}}
\end{equation} 
over the set $\II$. By the same reasoning in Theorem \ref{thm:mixedstrat} we have that there exists a fixed point of the set-valued map $E$ and by the symmetry of the game we deduce that the n-tuple $(\overline\mu, \dots, \overline\mu)$ is a Nash equilibrium in mixed strategies. \qed

\medskip

We call a Nash equilibrium $(\overline\mu_{j})_{j=1}^{N}$ in mixed strategies with $\overline\mu_{j}=\overline\mu_{k}$ for any $j, k \in \{1, \dots, N\}$ a {\it symmetric Nash equilibrium in mixed strategies}.

Hereafter, we consider a function $F$ and a coupling function $V$ of the form
\begin{equation*}
V^{i}[m_{1}, \dots, m_{N}](x) = \int_{\T^{d(N-1)}}{F\left(x_{i}, \frac{1}{N-1}\sum_{j \not= i} \delta_{x_{j}} \right)\ \prod_{\substack{j=1 \\ j \not= i}}^{N}m_{j}({\rm dx_{j}})}. 
\end{equation*}

\begin{proposition}\label{prop:measconv}
Let $(\overline\mu_{j}^{N})_{j=1}^{N} \in \II$ be  a Nash equilibrium in mixed strategies for the symmetric game, i.e., $\overline\mu^{N}_{i} = \overline\mu^{N}_{j}$ for any $j$, $i \in \{1, \dots, N\}$. Then, $\{\overline\mu^{N}\}_{N \in \N}$ is tight and there exists a probability measure $\overline\mu$ such that $\mu^{N} \rightharpoonup \overline\mu$. Moreover, we have that $V^{i}[\pi\sharp\overline\mu^{N}, \dots, \pi\sharp\overline\mu^{N}](x) \to F(x, \pi\sharp\overline\mu)$ as $N \to \infty$.  
\end{proposition}
\proof
Since $\overline\mu^{N} \in \II$, from the invariance w.r.t. the Euler flow and \cite[Corollary 4.4.5]{bib:FA} we easily get that $\{\overline\mu^{N}\}_{N \in \N}$ is tight. So, by Prokhorov's Theorem there exists a probability measure $\overline\mu$ such that $\mu^{N} \rightharpoonup \overline\mu$. Finally, by Theorem \ref{thm:HST} we obtain 
\begin{equation*}
\lim_{N \to \infty} V^{i}[\pi\sharp\overline\mu^{N}, \dots, \pi\sharp\overline\mu^{N}](x) = F(x, \pi \sharp \overline\mu). \eqno\square
\end{equation*}

For any $N \in \N$, from weak KAM theory we know that there exists a viscosity solution $v^{N}$ to the ergodic Hamilton-Jacobi equation 
\begin{equation}\label{eq:VHJ}
\lambda^{N} + H(x, Dv^{N}(x))=V^{N}[\pi\sharp\overline\mu^{N}, \dots, \pi\sharp\overline\mu^{N}](x), \quad x \in \T^{d}
\end{equation}
where $H$ is the Hamiltonian associated with $L$.

\begin{theorem}\label{thm:main}
For any $N \in \N$, let $(\lambda^{N}, v^{N}, \overline\mu^{N})$ be such that $(\lambda^{N}, v^{N})$ solves \eqref{eq:VHJ} and $\overline\mu^{N}$ is a symmetric Nash equilibrium in mixed strategies. Then, there exists $(\overline\lambda, \overline{v}, \overline\mu)$ such that
\begin{itemize}
\item[($i$)] $\overline\mu^{N}$ converges to $\overline\mu$;
\item[($ii$)] $\lambda^{N}$ converges to $\overline\lambda$ and we have that 
\begin{align*}
\overline\lambda =\ \inf_{\mu \in \II} \int_{\T^{d} \times \R^{d}}{\big(L(x,v) + F(x, \pi\sharp\overline\mu) \big)\ \mu({\rm dxdv})} =\  \int_{\T^{d} \times \R^{d}}{\big(L(x,v) + F(x, \pi\sharp\overline\mu) \big)\ \overline\mu({\rm dxdv})}.
\end{align*}
\item[($iii$)] $v^{N}$ uniformly converges to $\overline{v}$.
\end{itemize}
Moreover, $(\overline\lambda, \overline{v}, \overline\mu)$ solves the MFG system 
\begin{equation}\label{eq:ergoMFG}
\begin{cases}
\overline\lambda + H(x, D\overline{v}(x))=F(x, \pi\sharp\overline\mu), & x \in \T^{d}
\\
\ddiv\big(\pi\sharp\overline\mu D_{p}H(x, D\overline{v}(x)) \big)=0, & x \in \T^{d}.
\end{cases}
\end{equation}
\end{theorem}
\proof
Note that, ($i$) follows by Proposition \ref{prop:measconv}. Next, we now proceed to prove ($ii$). Since $\overline\mu^{N}$ converges to $\overline\mu$, it is enough to prove that $\overline\mu$ is a minimizer of the problem
\begin{equation*}
\inf_{\mu \in \II} \int_{\T^{d} \times \R^{d}}{\big(L(x,v) + F(x, \pi\sharp\overline\mu \big)\ \mu({\rm dxdv})}. 
\end{equation*}
The proof of this point is a simple adaptation of the one in Theorem \ref{thm:mixedstrat} and for this is reason we omit it here. 
In conclusion, from classical optimal control and weak KAM theory we have that $v^{N}$ uniformly converges to $\overline{v}$ and from the stability of viscosity solution  we deduce that 
\begin{equation*}
\overline\lambda + H(x, D\overline{v}(x))=F(x, \pi\sharp\overline\mu) \quad (x \in \T^{d}). \eqno\square
\end{equation*}

\section{Future direction}
\label{sec:Open}

From Theorem \ref{thm:main} we have that as the number of players $N$ goes to infinity the Nash equilibria in mixed strategies converge to a solution of a PDEs system defined by a continuum of agents. Heuristically speaking, it means that a solution to \eqref{eq:ergoMFG} is almost a Nash equilibrium for the $N$-player game. However, such a solution is not unique and this leads to a selection problem, i.e., we need some information more on the limit given in Theorem \ref{thm:main} besides being a solution to the system.

A possible approach toward a solution to such a problem can be inspired by \cite{Davini_2016} and \cite{A_Gomes_2020}, where the selection problem have been studied in case of Tonelli Hamiltonian system in the first work and in case of first-order mean field games in the second one. 

A different problem that can be address in a similar manner of this paper is the case of time dependent mean field game system. In a recent paper, \cite{bib:PC1} evolutive minimizing measures have been introduced and one can ask whether this flow of measures define a Nash equilibria in mixed strategies for the time dependent differential N-players game.

\appendix
 \section{Hewitt and Savage Theorem}

\begin{theorem}\label{thm:HST}
Let $\{m_{N}\}_{N \in \N}$ be a sequence of symmetric probability measures on $\T^{dN}$ such that 
\begin{equation*}
\int_{\T^{d}}{m_{N+1}({\rm dx})} = m_{N}
\end{equation*}
for all $n \in \N$. Then, there is a probability measure $\mu \in \PP(\T^{d})$ such that, for any test function $\phi \in C(\PP(\T^{d}))$ we have that 
\begin{equation*}
\lim_{N \to \infty} \int_{\T^{dN}}{\phi\left(\frac{1}{N}\sum_{i=1}^{N}{\delta_{x_{i}}}\right)\ m_{N}({\rm dx_{1}, \dots, dx_{N}})} = \int_{\PP(\T^{d})}{\phi(\mu)\ \mu({\rm dm})}.
\end{equation*}
\end{theorem}

%

\end{document}